\documentclass[12pt]{article}
 \usepackage[paperheight=12in,paperwidth=8.75in]{geometry}

\usepackage{setspace}
\usepackage[T1]{fontenc}
\usepackage{times}
\usepackage{palatino}
\usepackage{lmodern}
\usepackage[latin1]{inputenc}
\usepackage{epsfig}
\usepackage[english]{babel}
\usepackage{color}
\usepackage{graphicx}
\usepackage{dcolumn}
\usepackage{moreverb}
\usepackage{amsmath,amssymb,amsfonts}
\usepackage[all]{xy}

\begin{document}
\numberwithin{equation}{section}
\newcommand{\boxedeqn}[1]{%
  \[\fbox{%
      \addtolength{\linewidth}{-2\fboxsep}%
      \addtolength{\linewidth}{-2\fboxrule}%
      \begin{minipage}{\linewidth}%
      \begin{equation}#1\end{equation}%
      \end{minipage}%
    }\]%
}


\newsavebox{\fmbox}
\newenvironment{fmpage}[1]
     {\begin{lrbox}{\fmbox}\begin{minipage}{#1}}
     {\end{minipage}\end{lrbox}\fbox{\usebox{\fmbox}}}

\raggedbottom
\onecolumn

\newtheorem{th1}{Theorem}[section]
\newtheorem{df}{Definition}[section]
\newtheorem{lem}{Lemma}[section]
\newtheorem{cor}{Corollary}[section]
\newtheorem{pro}{Proposition}[section]
\newtheorem{ex}{Example}[section]
\newtheorem{rem}{Remark}[section]

\parindent 8pt
\parskip 10pt
\baselineskip 16pt
\noindent\title*{{\Large{\textbf{Left centralizers on Lie ideals in prime and semiprime gamma rings }}}}
\newline
\newline
Md Fazlul Hoque$^a$ and Akhil Chandra Paul$^b$
\newline
\newline
$a.$ Department of Mathematics, Pabna University of Science and Technology, Pabna-6600, Bangladesh
\newline
\newline
$b.$ Department of Mathematics, University of Rajshahi, Rajshahi-6205, Bangladesh
\newline
\newline
E-mail: fazlul\_math@yahoo.co.in; acpaulrubd\_math@yahoo.com
\newline
\newline

\begin{abstract}
Let $U$ be a Lie ideal of a 2-torsion free prime $\Gamma $-ring $M$ such that $u\alpha u\in U$ for all $u\in U$ and $\alpha \in\Gamma $. 
If $T:M\rightarrow M$ is an additive mapping satifying the relation $T(u\alpha u)=T(u)\alpha u$ for all $u\in U$ and $\alpha \in\Gamma $, then we prove that $T(u\alpha v)=T(u)\alpha v$ for all $u, v\in U$ and $\alpha \in\Gamma $.
Also this result is extended to semiprime $\Gamma $-rings.  . 
\end{abstract}



\section{Introduction}
The concept of a gamma ring was presented as a generalization of the classical rings. This gamma ring was first introduced by Nobusawa\cite{Nob1} in 1964,  which is currently known as $\Gamma_N$-ring. After two years it was more broadly generalized by Bernes \cite{Ber1} in the sense of Nobusawa\cite{Nob1}, which is now known as the $\Gamma$-ring. It is shown that a $\Gamma$-ring need not be a ring, but a $\Gamma$-ring is more general than rings \cite{Nob1}, and also that every $\Gamma_N$-ring is a $\Gamma$-ring \cite{Ber1}. From its beginning, the various important theories of the classical rings were extended and generalized to the theories of $\Gamma$-rings \cite{Luh1, Kyu1}. Such theories have been attracted much international attentions as an emerging areas of research to the modern algebraists to enrich the areas of algebras. Recently, many researchers determine a number of basic properties of $\Gamma$-rings with creative and productive remarkable results \cite{Cev1, Sab1, Hoq1, Hal1, Ull1, Dey1, Hoq4}.

Borut Zalar \cite{Zal1} studied on centralizers of semiprime rings and shown that Jordan centralizers and centralizers of this rings coincide. Using the concept of centralizers, Vukman \cite{Vuk1,Vuk2} established a number of results on prime and semiprime rings. Such results and ideas have been extended to prime and semiprime $\Gamma$-rings in different aspects such as centralizers and $\theta$-centralizers \cite{Hoq1, Ull1, Hoq4, Hoq2, Hoq5}, Jordan centralizers and Jordan $\theta$-centralizers \cite{Hoq3, Sal1} and centralizers with involutions \cite{Hoq6}.

The Lie ideals and Jordan derivations of prime rings was studied in \cite{Awt1, Ash1, Cor1, Hon1}, and these works have been extended to Lie ideals of prime $\Gamma$-rings of Jordan derivations \cite{Hal1, Hon1, Reh1} and Jordan $k$-derivations \cite{Naz1}. In fact, a number of significant results of classical ring theories were developed in prime and semiprime $\Gamma$-rings with Lie ideals and Jordan structures \cite{Sab1, Rah1}. However, the research on centralizers of prime and semiprime gamma rings with Lie ideals is still an unexplored area and it would be of interest to further works. Thus the aim of this article is to extend the results of \cite{Hoq1} to the Lie ideals in prime and semiprime $\Gamma $-rings.

\section{Preliminaries}
Let $M$ and $\Gamma $  be additive abelian groups. If there exists a mapping $(a,\alpha ,b)\rightarrow  a\alpha b$ of $M\times\Gamma \times M\rightarrow M$, which satisfies the conditions, 
\begin{itemize}
\setlength{\itemsep}{0pt}
  \setlength{\parskip}{0pt}
  \setlength{\parsep}{0pt}
\item[(i)]
$a\alpha b\in M$
\item[(ii)] $(a+b)\alpha c$=$a\alpha c$+$b\alpha c$, $a(\alpha +\beta )c$=$a\alpha c$+$a\beta c$, $a\alpha (b+c)$=$a\alpha b$+$a\alpha c$
\item[(iii)]
$(a\alpha b\beta c$=$a\alpha (b\beta c)$ for all $a,b,c\in M$ and $\alpha ,\beta \in \Gamma$,
\end{itemize}
then $M$ is called a $\Gamma$-ring.
\\Let $M$ be a $\Gamma $-ring. Then $M$ is said to be prime if $a\Gamma M\Gamma b$=$(0)$ with $a,b\in M$, implies $a$=$0$ or $b$=$0$, and semiprime if $a\Gamma M\Gamma a$=$(0)$ with $a\in M$ implies $a$=$0$. An additive subgroup $U$ of $M$ is said to be a Lie ideal of $M$ if $[u,x]_{\alpha }\in U$ for all $u\in U$, $x\in M$ and $\alpha \in\Gamma $.
Furthermore, $M$ is said to be commutative $\Gamma $-ring if $a\alpha b$=$b\alpha a$ for all $a, b\in M$ and $\alpha \in\Gamma $. 
Moreover, the set $Z(M)$ =$\{a\in M:a\alpha b=b\alpha a $ for all $ \alpha \in \Gamma,  b\in M\}$ is called the centre of the $\Gamma $-ring $M$.
\\If $M$ is a $\Gamma $-ring, then $[a,b]_{\alpha }$=$a\alpha b-b\alpha a$ is known as the commutator of $a$ and $b$ with respect to $\alpha $, where $a, b\in M$ and $\alpha \in\Gamma $. It has the basic commutator identities:
\begin{eqnarray*}
 &&[a\alpha b,c]_{\beta }=[a,c]_{\beta }\alpha b+a[\alpha ,\beta ]_{c}b+a\alpha [b,c]_{\beta },
 \\&&[a,b\alpha c]_{\beta }=[a,b]_{\beta }\alpha c+b[\alpha ,\beta ]_{a}c+b\alpha [a,c]_{\beta },
\end{eqnarray*}
for all $a,b, c\in M$ and $\alpha ,\beta \in\Gamma $.
One consider the following assumption \cite{Hoq1},
\begin{eqnarray*}
(A)\quad.................\quad\alpha b\beta c=a\beta b\alpha c,
\end{eqnarray*}
for all $a,b,c\in M$, and $\alpha ,\beta \in\Gamma $, which will extensively use through the paper. According to the assumption $(A)$, the above two identites reduce to
\begin{eqnarray*}
 &&[a\alpha b,c]_{\beta }=[a,c]_{\beta }\alpha b+a\alpha [b,c]_{\beta },
 \\&&[a,b\alpha c]_{\beta }=[a,b]_{\beta }\alpha c+b\alpha [a,c]_{\beta }.
\end{eqnarray*}
For existence of such a $\Gamma $-ring $M$, we present the following example.
\begin{ex}(\cite{Cev1}, Example 1.1) Let $R$ be an associative ring with the unity element 1. Let $M=M_{1, 2}(R)$ and $\Gamma =\left\{\left(\begin{array}{c}n.1\\0\end{array}\right): \mbox{n is an integer}\right\}$. Then $M$ is a $\Gamma $-ring.
A simple verification shows that $M$ satisfies the assumption (A).
\end{ex} 
An additive mapping $T:M\rightarrow M$ is called a left (right) centralizer if $T(a\alpha b)=T(a)\alpha b$ (resp. $T(a\alpha b)=a\alpha T(b)$) for all $a, b\in M$ and $\alpha \in \Gamma $.
 A centralizer is an additive mapping which is both a left and a right centralizer. For any fixed $a\in M$ and $\alpha \in\Gamma$, the mapping $T(x)=a\alpha x$ is a left centralizer, and $T(x)=x\alpha a$ is a right centralizer. We shall restrict our attention on left centralizer, since all results of right centralizers are the same as left centralizers.
An additive mapping $T:M\rightarrow M$ is Jordan left(right) centralizer if $T(x\alpha x)=T(x)\alpha x(T(x\alpha x)=x\alpha T(x))$ for all $x\in M$ and $\alpha \in\Gamma $.
Every left centralizer is a Jordan left centralizer but the converse is not in general true.
An additive mappings $T:M\rightarrow M$ is called a Jordan centralizer if $T(x\alpha y+y\alpha x)=T(x)\alpha y+y\alpha T(x)$, for all $x,y\in M$ and $\alpha \in\Gamma $. Every centralizer is a Jordan centralizer but Jordan centralizer is not in general a centralizer.

\section{Left centralizers of prime gamma rings}
\setlength{\itemsep}{0pt}
  \setlength{\parskip}{0pt}
  \setlength{\parsep}{0pt}
\begin{lem}\label{lem1}
Let $M$ be a $\Gamma$-ring and $U$ a Lie ideal of $M$ such that $u\alpha u\in U$ for all $u\in U$ and $\alpha \in\Gamma $.
If $T:M\rightarrow M$ is an additive mapping satisfying the relation $T(u\alpha u)=T(u)\alpha u$ for all $u\in U$ and $\alpha \in\Gamma $, then
\begin{itemize}
\setlength{\itemsep}{0pt}
  \setlength{\parskip}{0pt}
  \setlength{\parsep}{0pt}
\item[(a)]
$T(u\alpha v+v\alpha u)=T(u)\alpha v+T(v)\alpha u;$
\item[(b)]
$T(u\alpha v\beta u+u\beta v\alpha u)=T(u)\alpha v\beta u+T(u)\beta v\alpha u;$
\item[(c)]
$T(u\alpha v\beta u)=T(u)\alpha v\beta u$; 
\item[(d)]
$T(u\alpha v\beta w+w\beta v\alpha u)=T(u)\alpha v\beta w+T(w)\beta v\alpha u$,
\end{itemize}
 for all $u, v, w\in U$ and $\alpha ,\beta \in\Gamma $.      
\end{lem}
{\bf Proof.} By the definition of Lie ideal $U$, $u\alpha u\in U$ for all $u\in U$ and $\alpha \in\Gamma $. Thus we have,  
$u\beta v+v\beta u=(u+v)\beta (u+v)-u\beta u-vv\in U$ for all $u, v\in U$ and $\beta \in\Gamma $. 
Therefore \begin{eqnarray*}T(u\alpha v+v\alpha u)&&=T((u+v)\alpha (u+v))-T(u\alpha u)-T(v\alpha v)\\&&=T(u+v)\alpha (u+v)-T(u)\alpha u-T(v)\alpha v
\\&&=T(u)\alpha u+T(u)\alpha v+T(v)\alpha u+T(v)\alpha v-T(u)\alpha -T(v)\alpha v\\&&=T(u)\alpha v+T(v)\alpha u.\end{eqnarray*} 
Hence 
\begin{eqnarray}
&T(u\alpha v+v\alpha u)&= T(u)\alpha v+T(v)\alpha u.\label{eq1}
\end{eqnarray} 
Since $u\beta v+v\beta u\in U$ for all $u, v\in U$ and $\beta \in\Gamma $, we replace $v$ by $u\beta v+v\beta u$ in relation (\ref{eq1}), we obtain
\begin{eqnarray*}T(u\alpha (u\beta v+v\beta u)+(u\beta v+v\beta u)\alpha u)&=&T(u)\alpha (u\beta v+v\beta u)+T(u\beta v+v\beta u)\alpha u\end{eqnarray*}
\begin{eqnarray*}\lefteqn{\Rightarrow T(u\alpha u\beta v+u\alpha v\beta u+u\beta v\alpha u+v\beta u\alpha u)}\\&&=T(u)\alpha (u\beta v+v\beta u)+T(u)\beta v\alpha u+T(v)\beta u\alpha u.\end{eqnarray*}
\begin{eqnarray*}\lefteqn{\Rightarrow T(u\alpha v\beta u+u\beta v\alpha u)+T(u\alpha u)\beta v+T(v)\beta u\alpha u}\\&&=T(u)\alpha u\beta v+T(u)\alpha v\beta u+T(u)\beta v\alpha u+T(v)\beta u\alpha u\end{eqnarray*}     
\begin{eqnarray}
\Rightarrow T(u\alpha v\beta u+u\beta v\alpha u)&=&T(u)\alpha v\beta u+T(u)\beta v\alpha u.\label{eq2}\end{eqnarray} 
Hence (b) proved.
\\By using the assumption $(A)$ in the above relation (\ref{eq2}), we obtain 
 \begin{eqnarray*}
 2T(u\alpha v\beta u)&=&2T(u)\alpha v\beta u.
 \end{eqnarray*} 
Thus for the 2-torsion freeness of $M$, we have
 \begin{eqnarray}
 T(u\alpha v\beta u)&=&T(u)\alpha v\beta u. \label{eq3}
 \end{eqnarray} 
Putting $u=u+w$ in the relation (\ref{eq3}), we obtain the result (d).
\\\\Define $B_{\alpha }(u, v)=T(u\alpha v)-T(u)\alpha v$ for all $u, v\in U$ and $\alpha \in\Gamma $. Then we have the following remarks and lemmas.
\begin{rem} It is clear that $B_{\alpha }(u, v)$ is an additive mapping such that $B_{\alpha }(u, v)+B_{\alpha }(v, u)=0$. 
\end{rem} 
\begin{rem} It is also clear that $T$ is a left centralizer if and only if $B_{\alpha }(u, v)=0$. 
\end{rem} 
\begin{lem}\label{lem2}
Let $M$ be a 2-torsion free $\Gamma $-ring and $U$ be a Lie ideal of $M$ such that $u\alpha u\in U$ for all $u\in U$ and $\alpha \in\Gamma $.
If $T:M\rightarrow M$ is an additive mapping satisfying the relation $T(u\alpha u)=T(u)\alpha u$ for all $u\in U$ and $\alpha \in\Gamma $, then 
$B_{\alpha }(u, v)\beta w\gamma [u, v]_{\delta }=0$ and $[u, v]_{\delta }\beta w\gamma B_{\alpha }(u, v)=0$. 
\end{lem}    
{\bf Proof.} First we compute 
\begin{eqnarray}
x&=&T(2(u\alpha v)\beta w\gamma 2(v\delta u)+2(v\alpha u)\beta w\gamma 2(u\delta v))\label{eqx1}
\end{eqnarray} 
in two different ways. Then using Lemma-\ref{lem1} (c) in (\ref{eqx1}), we have
\begin{eqnarray}
x&=&4T(u)\alpha v\beta w\gamma v\delta u+4T(v)\alpha u\beta w\gamma u\delta v,\label{eq4}
\end{eqnarray}   
and using Lemma-\ref{lem1} (d) in (\ref{eqx1}), we have
\begin{eqnarray}
x&=&4T(u\alpha v)\beta w\gamma v\delta u+4T(v\alpha u)\beta w\gamma u\delta v\label{eq5}
\end{eqnarray}
Comparing (\ref{eq4}) and (\ref{eq5}), we obtain 
\begin{eqnarray*}
4\{(T(u\alpha v)-T(u)\alpha v)\beta w\gamma v\delta u+(T(v\alpha u)-T(v)\alpha u)\beta w\gamma u\delta v\}&=&0.
\end{eqnarray*}
\begin{eqnarray*}
\Rightarrow 4\{B_{\alpha }(u, v)\beta w\gamma v\delta u+B_{\alpha }(v, u)\beta w\gamma u\delta v\}&=&0
\end{eqnarray*}
\begin{eqnarray*}
\Rightarrow 4\{B_{\alpha }(u, v)\beta w\gamma v\delta u-B_{\alpha }(u, v)\beta w\gamma u\delta v\}&=&0.
\end{eqnarray*}
Hence we have 
\begin{eqnarray*}
4B_{\alpha }(u, v)\beta w\gamma [u, v]_{\delta}&=&0.
\end{eqnarray*}
Therefore, by the semiprimeness of $M$, we obtain 
\begin{eqnarray*}
B_{\alpha }(u, v)\beta w\gamma [u, v]_{\delta}&=&0
\end{eqnarray*} 
for all $u, v, w\in U$ and $\alpha ,\beta ,\gamma ,\delta \in\Gamma $. 
\\Similarly, we can easily prove that 
\begin{eqnarray*}[u, v]_{\delta}\beta w\gamma B_{\alpha }(u, v)&=&0
\end{eqnarray*} for all $u, v, w\in U$ and $\alpha ,\beta ,\gamma ,\delta \in\Gamma $.  
\begin{lem}(\cite{Hal1},Lemaa-2)\label{lem3}
 Let $U\not\subseteq Z(M)$ be a Lie ideal of a 2-torsion free prime $\Gamma $-ring $M$ and $a, b\in M$ such that $a\alpha U\beta b=0$. Then $a=0$ or $b=0$.  
\end{lem}  
\begin{lem} \label{lem4}
Let $U$ be a commutative Lie ideal of a 2-torsion free prime $\Gamma $-ring $M$. Then $U\subseteq Z(M)$.  
\end{lem}
{\bf Proof.} For $u\in U$, $m\in M$ and $\alpha \in\Gamma $, we have $[u, m]_{\alpha }\in U$. Since $U$ is commutative, $[u, [u, m]_{\alpha }]_{\beta }=0$ for all $\beta \in\Gamma $.
Now for $x, y\in M$ and $\gamma \in\Gamma $, replace $x\gamma y$ for $m$, we obtain $[u, [u, x\gamma y]_{\alpha }]_{\beta }=0$. By using $(A)$, we have $[u, x\gamma [u, y]_{\alpha }+[u, x]_{\alpha }\gamma y]_{\beta }=0$, 
$\Rightarrow [u, x\gamma [u, y]_{\alpha }]_{\beta }+[u, [u, x]_{\alpha }\gamma y]_{\beta }=0$,
$\Rightarrow x\gamma [u, [u, y]_{\alpha }]_{\beta }+[u, x]_{\beta }\gamma [u, y]_{\alpha }+[u, x]_{\alpha }\gamma [u, y]_{\beta }+[u, [u, x]_{\alpha }]_{\beta }\gamma y=0$,
$\Rightarrow [u, x]_{\beta }\gamma [u, y]_{\alpha }+[u, x]_{\alpha }\gamma [u, y]_{\beta }=0$. After some calculation and using the assumption $(A)$, we have
$2[u, x]_{\alpha }\gamma [u, y]_{\beta }=0$. Since $M$ is 2-torsion free, thus $[u, x]_{\alpha }\gamma [u, y]_{\beta }=0$. 
\\Putting $y$ by $y\delta m$ for all $m\in M$, we have $[u, x]_{\alpha }\gamma [u, y\delta m]_{\beta }=0$. This implies $[u, x]_{\alpha }\gamma y\delta [u, m]_{\beta }+[u, x]_{\alpha }\gamma [u, y]_{\beta }\delta m=0$, 
 $\Rightarrow [u, x]_{\alpha }\gamma y\delta [u, m]_{\beta }=0$. Hence by primeness of $M$, $[u, x]_{\alpha }=0$ or $[u, m]_{\beta }=0$. 
If $[u, x]_{\alpha }=0$, then $U\subseteq Z(M)$ and if $[u, m]_{\beta }=0$, then also $U\subseteq Z(M)$. 

\vskip.2cm
\noindent 
\begin{th1} \label{th1}
Let $U$ be a Lie ideal of a 2-torsion free prime $\Gamma $-ring $M$ such that $u\alpha u\in U$ for all $u\in U$ and $\alpha \in\Gamma $.
If $T:M\rightarrow M$ is an additive mapping such that $T(u\alpha u)=T(u)\alpha u$ for all $u\in U$ and $\alpha \in\Gamma $, then $T(u\alpha v)=T(u)\alpha v$ for all $u, v\in U$ and $\alpha \in\Gamma $.
\end{th1}
\vskip.2cm
\noindent   
{\bf Proof.} If $U$ is a commutative Lie ideal of $M$, then by Lemma-\ref{lem4}, $U\subseteq Z(M)$. Therefore, by Lemma-\ref{lem1} (d), 
we have 
\begin{eqnarray}
&T(u\alpha v\beta w+w\beta v\alpha u)&=T(u)\alpha v\beta w+T(w)\beta v\alpha u.\label{eq6}
\end{eqnarray}  
Since $U$ is commutative, we have $u\alpha v=v\alpha u$. Therefore 
\begin{eqnarray}
&T((u\alpha v)\beta w+w\beta (u\alpha v))&=T(u\alpha v)\beta w+T(w)\beta u\alpha v.\label{eq7}
\end{eqnarray} 
Comparing (\ref{eq6}) and (\ref{eq7}), and using $u\alpha v=v\alpha u$, 
we obtain 
\begin{eqnarray*}
&(T(u\alpha v)-T(u)\alpha v)\beta w&=0,
\end{eqnarray*} 
which yields that $G_{\alpha }(u,v)\beta w=0$. 
Since $w\in U$, $[w,m]_{\gamma }\in U$ for all $m\in M$ and $\gamma \in\Gamma $. 
Replacing $w$ by $[w,m]_{\gamma }$, we obtain  $G_{\alpha }(u,v)\beta [w,m]_{\gamma }=0$. This implies $G_{\alpha }(u,v)\beta w\gamma m-G_{\alpha }(u,v)\beta m\gamma w=0$.
Hence the relation becomes $G_{\alpha }(u,v)\beta m\gamma w=0$ for all $u, v, w\in U$, $m\in M$ and $\alpha ,\beta ,\gamma \in\Gamma $. Since $U\not=0$, in view of the Lemma-\ref{lem3}, $G_{\alpha }(u,v)=0$.   
\\If $U$ is not commutative, then $U\not\subseteq Z(M)$. In this case, we have from Lemma-\ref{lem2}, $B_{\alpha }(u, v)\beta w\gamma [u, v]_{\delta }=0$.
Putting $u=u+x$ for all $u\in U$, we have $(B_{\alpha }(u, v)+B_{\alpha }(x, v))\beta w\gamma ([u, v]_{\delta }+[x, v]_{\delta })=0$.
This implies $B_{\alpha }(u, v)\beta w\gamma [x, v]_{\delta }+B_{\alpha }(x, v)\beta w\gamma [u, v]_{\delta }=0$. 
\\Now, $B_{\alpha }(u, v)\beta w\gamma [x, v]_{\delta }\mu z\gamma B_{\alpha }(u, v)\beta w\gamma [x, v]_{\delta }$
\\$=-B_{\alpha }(u, v)w\gamma [x, v]_{\delta }\mu z\gamma B_{\alpha }(x, v)\beta w\gamma [u, v]_{\delta }$
\\$=0$. 
\\Therefore, by Lemma-\ref{lem3}, we have $B_{\alpha }(u, v)\beta w\gamma [x, v]_{\delta }=0$ for all $x\in U$. Similarly, using $v=v+y$, we obtain $B_{\alpha }(u, v)\beta w\gamma [x, y]_{\delta }=0$ for all $y\in U$.
By Lemma-\ref{lem3}, we obtain $B_{\alpha }(u, v)=0$ or $[x, y]_{\delta }=0$. If $[x, y]_{\delta }=0$, then $U$ is commutative which shows a contradiction that $U\not\subseteq Z(M)$.
Therefore $B_{\alpha }(u, v)=0$.  
\begin{cor} Let $M$ be a 2-torsion free prime $\Gamma $-ring and $T:M\rightarrow M$ be a Jordan left centralizer. Then $T$ is a left centralizer.
\end{cor}
\vskip.2cm   
\noindent

\section{Left centralizers of semiprime gamma rings}  
\begin{lem} \label{lem41}
Let $U$ be a commutative Lie ideal of a 2-torsion free semiprime $\Gamma $-ring $M$. Then $U\subseteq Z(M)$.   
\end{lem}          
{\bf Proof.} For $u\in U$ and $x\in M$, we have $[u, [u, x]_{\alpha }]_{\alpha }=0$. Repacing $x=x\gamma y$, we have $[u, [u, x\gamma y]_{\alpha }]_{\alpha }=0$.
This implies $[u, x\gamma [u, y]_{\alpha }+[u, x]_{\alpha }\gamma y]_{\alpha }=0$,
 $\Rightarrow [u, x]_{\alpha }\gamma [u, y]_{\alpha }+x\gamma [u, [u,y]_{\alpha }]_{\alpha }+[u, [u,x]_{\alpha }]_{\alpha }\gamma y+[u,x]_{\alpha }\gamma [u,y]_{\alpha }=0$,
$\Rightarrow [u, x]_{\alpha }\gamma [u, y]_{\alpha }+[u,x]_{\alpha }\gamma [u,y]_{\alpha }=0$, i.e., $2[u, x]_{\alpha }\gamma [u, y]_{\alpha }=0$. Hence by 2-torsion freeness, we have $[u, x]_{\alpha }\gamma [u, y]_{\alpha }=0$.
Replacing $y$ by $y\delta x$, we have $[u, x]_{\alpha }\gamma [u, y\delta x]_{\alpha }=0$. This implies, $[u, x]_{\alpha }\gamma [u, y]_{\alpha }\delta x+[u, x]_{\alpha }\gamma y\delta [u, x]_{\alpha }=0$, that is, 
$[u, x]_{\alpha }\gamma y\delta [u, x]_{\alpha }=0$ for all $y\in M$. Since $M$ is semiprime, $[u,x]_{\alpha }=0$, which shows that $U\in Z(M)$.        
\begin{lem} \label{lem42}
Let $U$ be a Lie ideal of a 2-torsion free $\Gamma $-ring $M$ satisfying the assumption $(A)$, then $T(U)=\{x\in M : [x,M]_{\Gamma }\subseteq U\}$ is both a subring
and a Lie ideal of $M$ such that $U\subseteq T(M)$. 
\end{lem}
{\bf Proof.} Since $U$ is a Lie ideal of $M$, so we have $[U,M]_{\Gamma }\subseteq U$. Thus $U\subseteq T(U)$. Also, we have $[T(U),M]_{\Gamma }\subseteq U\subseteq T(U)$. Hence $T(U)$ is a Lie ideal of $M$. 
\\Suppose that $x, y\in T(U)$, then $[x,m]_{\alpha }\in U$ and $[y,m]_{\alpha }\in U$, for all $m\in M$ and $\alpha \in\Gamma $.  
\\Now, $[x\alpha y,m]_{\beta }=x\alpha [y,m]_{\beta }+[x,m]_{\beta }\alpha y\in U$. Hence $[x\alpha y,m]_{\beta }\in U$, for all $x, y\in T(U)$, $m\in M$ and $\alpha ,\beta \in\Gamma $. Therefore, $x\alpha y\in T(U)$.
\begin{lem} \label{lem43}
Let $U\not\subseteq Z(M)$ be a Lie ideal of a 2-torsion free semiprime $\Gamma $-ring $M$. Then there exists a nonzero ideal $K=M\Gamma [U,U]_{\Gamma }\Gamma M$ of $M$ generated by $[U,U]_{\Gamma }$ such that $[K,M]_{\Gamma }\subseteq U$.  
\end{lem}
{\bf Proof.} First, we have to prove that if $[U,U]_{\Gamma }=0$, then $U\subseteq Z(M)$. Let $[U,U]_{\Gamma }=0$. Then for all $a\in U$, $\alpha \in\Gamma $,
we have $[u,[u,x]_{\alpha }]_{\alpha }=0$ for all $x\in M$. Then using the proof of Lemma-\ref{lem41}, we obtain $U\subseteq Z(M)$, which is a contradiction. Thus, let $[U,U]_{\Gamma}\not=0$.
Then $K=M\Gamma [U,U]_{\Gamma }\Gamma M$ is a nonzero ideal of $M$ generated by $[U,U]_{\Gamma }$. Let $x, y\in U$, $m\in M$ and $\alpha , \beta \in\Gamma $, 
we have $[x,y\beta m]_{\alpha }, y, [x,m]_{\alpha }\in U\subseteq T(U)$. Hence by Lemma-\ref{lem42}, $[x, y]_{\alpha }\beta m=[x,y\beta m]_{\alpha }-y\beta [x,m]_{\alpha }\in T(U)$. 
\\Also, we can show that $m\beta [x,y]_{\alpha } \in T(U)$ and therefore, we obtain  $[[U,U]_{\Gamma },M]_{\Gamma }\subseteq U$. 
That is, $[[[[x,y]_{\alpha },m]_{\alpha },s]_{\alpha },t]_{\alpha }\in U$ for all $m, s, t\in M$ and $\alpha \in\Gamma $. 
Hence $[[x,y]_{\alpha }\alpha m\alpha s-m\alpha [x,y]_{\alpha }\alpha s+[s,m]_{\alpha }\alpha [x,y]_{\alpha }-[s\alpha [x,y]_{\alpha },m]_{\alpha },t]_{\alpha }\in T(U)$.    
Since $[x, y]_{\alpha }\alpha m\alpha s$, $s\alpha [x,y]_{\alpha }$, $[s,m]_{\alpha }\alpha [x,y]_{\alpha }\in T(U)$. Thus we have $[m\alpha [x,y]_{\alpha }\alpha s,t]_{\alpha }\in U$ for all $m, s, t\in M$ and $\alpha \in\Gamma $.
Hence $[K,M]_{\Gamma }\subseteq U$.
\begin{lem}\label{lem44}
 Let $U\not\subseteq Z(M)$ be a Lie ideal of a 2-torsion free semiprme $\Gamma $-ring $M$ and $a\in U$. If $a\alpha U\beta a=\{0\}$ for all $\alpha ,\beta \in\Gamma $,
then $a\alpha a=0$ and there exists a nonzero ideal $K=M\Gamma [U,U]_{\Gamma }\Gamma M$ of $M$ generated by $[U,U]_{\Gamma }$ such that $[K,M]_{\Gamma }\subseteq U$ and 
$K\Gamma a=a\Gamma K=\{0\}$.
\end{lem}  
{\bf Proof.} If $a\alpha U\beta a=\{0\}$ for all $\alpha ,\beta \in\Gamma $, then $a\alpha [a,a\delta m]_{\alpha }\beta a=0$ for all $m\in M$ and $\delta \in\Gamma $.
Therefore, \begin{eqnarray*}&0&=a\alpha (a\alpha a\delta m-a\delta m\alpha a)\beta a\\&&=a\alpha a\alpha a\delta m\beta a-a\alpha a\delta m\alpha a\beta a
\\&&=a\alpha a\delta a\alpha m\beta a-a\alpha a\delta m\beta a\alpha a\end{eqnarray*} Since $a\alpha a\delta a=0$, we have $(a\alpha a)\delta m\beta (a\alpha a)=0$ 
and hence $a\alpha a=0$ for semiprimeness of $M$. Now, we obtain $a\alpha [k\gamma a,m]_{\mu }\beta u\alpha a=0$ for all $k\in K$, $m\in M$, $u\in U$ and $\alpha ,\beta ,\mu \in\Gamma $. 
Therefore \begin{eqnarray*}&0&=a\alpha (k\gamma a\mu m-m\mu k\gamma a)\beta u\alpha a\\&&=a\alpha k\gamma a\mu m\beta u\alpha a-a\alpha m\mu k\gamma a\beta u\alpha a\\&&=a\alpha k\gamma a\mu m\beta u\alpha a\end{eqnarray*} 
Thus, we have $a\alpha k\gamma a\mu m\beta [k,a]_{\gamma }\alpha a=0$. This implies that $a\alpha k\gamma a\mu m\beta (k\gamma a-a\gamma k)\alpha a=0$ and 
hence $a\alpha k\gamma a\mu m\beta k\gamma a\alpha a-a\alpha k\gamma a\mu m\beta a\gamma k\alpha a=0$. Hence by using (A) and $a\alpha a=0$, we have $(a\alpha k\gamma a)\mu m\beta (a\alpha k\gamma a)=0$. 
Hence $a\alpha k\gamma a=0$ for semiprimeness of $M$. Thus we find that $(a\alpha k)\Gamma M\Gamma (a\alpha k)=0$. Hence $a\alpha k=0$ for all $k\in K$, that is, $a\alpha K=\{0\}$. Similarly, we have $K\alpha a=\{0\}$.           
 
\begin{lem}\label{lem45}
 Let $U\not\subseteq Z(M)$ be a Lie ideal of a 2-torsion free semiprme $\Gamma $-ring $M$ and $a, b\in U$ and $\alpha ,\beta \in\Gamma $.
\\{\bf (i)} If $a\alpha U\beta a=\{0\}$, then $a=0$.
\\{\bf (ii)} If $a\alpha U=\{0\}$ ($U\alpha a=\{0\}$), then $a=0$. 
\\{\bf (iii)} If $u\alpha u\in U$ for all $u\in U$ and $a\alpha U\beta b=\{0\}$, then $a\alpha b=0$ and $b\alpha a=0$ for all $\alpha \in\Gamma $.
\end{lem}
{\bf Proof.} (i) By Lemma-\ref{lem44}, we have $K\alpha a=M\Gamma [U,U]_{\Gamma }\Gamma M\alpha a=\{0\}$ and $a\alpha a=0$ for all $\alpha \in\Gamma $. Thus for all $x, y\in M$ and $\alpha,\beta \in\Gamma $, we have
\begin{eqnarray*}
&0&=[[a,x]_{\alpha },a]_{\gamma }\beta y\alpha a\\&&=[a\alpha x-x\alpha a,a]_{\gamma }\beta y\alpha a\\&&=a\alpha [x,a]_{\gamma }\beta y\alpha a-[x,a]_{\gamma }\alpha a\beta y\alpha a
\\&&=a\alpha x\gamma a\beta y\alpha a-a\alpha a\gamma x\beta y\alpha a-x\gamma a\alpha a\beta y\alpha a+a\gamma x\alpha a\beta y\alpha a\\&&=a\alpha x\gamma a\beta y\alpha a+a\gamma x\alpha a\beta y\alpha a\\&&=2a\alpha x\gamma a\beta y\alpha a.
\end{eqnarray*}  
By the 2-torsion freeness of $M$, we have, $a\alpha x\gamma a\beta y\alpha a=0$. Hence, we obtain $a\alpha x\gamma a\beta y\alpha a\delta x\gamma a=0$. By using the assumption (A), we have
$(a\alpha x\gamma a)\beta y\delta (a\alpha x\gamma a)=0$ for all $y\in M$. This implies $(a\alpha x\gamma a)\beta M\delta (a\alpha x\gamma a)=0$. Hence, by semiprimeness of $M$, we have $a\alpha x\gamma a=0$ for all $x\in M$ and $\alpha ,\gamma \in\Gamma $. Again, by semiprimeness of $M$, we obtain $a=0$.    
\\\\(ii) If $a\alpha U=\{0\}$, then $a\alpha U\beta a=\{0\}$ for all $\beta \in\Gamma $. Thus, by (i), we have $a=0$. Semilarly, if $U\alpha a=\{0\}$, then $a=0$. 
\\\\(iii) If $a\alpha U\beta b=\{0\}$, then we have $(b\gamma a)\alpha U\beta (b\gamma a)=\{0\}$ and hence, by (i), $b\gamma a=0$, for all $\gamma \in\Gamma $. 
Also, $(a\gamma b)\alpha U\beta (a\gamma b)=\{0\}$ if $a\alpha U\beta b=\{0\}$ and hence $a\gamma b=0$.  

\begin{th1} Let $U\not\subseteq  Z(M)$ be a Lie ideal of a 2-torsion free semiprime $\Gamma $-ring such that $u\alpha u\in U$ for all $u\in U$ and $\alpha \in\Gamma $. If $T:M\rightarrow M$ be an additive mapping satisfying the relation
$T(u\alpha u)=T(u)\alpha u$ for all $u\in U$ and $\alpha \in\Gamma $, then $T(u\alpha v)=T(u)\alpha v$ for all $u, v\in U$ and $\alpha \in\Gamma $. 
\end{th1}  
\vskip.2cm 
\noindent 
{\bf Proof.} Since $U\not\subseteq Z(M)$, we have from Lemma-\ref{lem2}, 
\begin{eqnarray*}
&B_{\alpha }(u,v)\beta w\gamma [u,v]_{\delta }&=0
\end{eqnarray*}   
By linearing $u$, we obtain 
\begin{eqnarray*}
&B_{\alpha }(u,v)\beta w\gamma [x,v]_{\delta }+B_{\alpha }(x,v)\beta w\gamma [u,v]_{\delta }&=0
\end{eqnarray*}
 for all $x\in U$. 
Now, \begin{eqnarray*}
\lefteqn{B_{\alpha }(u,v)\beta w\gamma [x,v]_{\delta }\mu z\nu B_{\alpha }(u,v)\beta w\gamma [x,v]_{\delta }}\\&&=-B_{\alpha }(u,v)\beta w\gamma [x,v]_{\delta }\mu z\nu B_{\alpha }(x,v)\beta w\gamma [u,v]_{\delta }\\&&=0
\end{eqnarray*} 
for all $z\in U$. Hence, by Lemma-\ref{lem45}(i), 
\begin{eqnarray*}
&B_{\alpha }(u,v)\beta w\gamma [x,v]_{\delta }&=0.
\end{eqnarray*} 
Similarly, linearing $v$, we obtain 
\begin{eqnarray*}
&B_{\alpha }(u,v)\beta w\gamma [x,y]_{\delta }&=0.
\end{eqnarray*} for all $y\in U$.        
Hence the similar proof of the Theorem-2.1 in \cite{Hoq1}, we obtain the required result.
\vskip.2cm 
\noindent

\begin{cor} 
Let $M$ be a 2-torsion free semiprime $\Gamma $-ring and $T:M\rightarrow M$ be a Jordan left centralizer. Then $T$ is a left centralizer.
\end{cor}
     
\begin{ex}
 Let $R$ be a commutative ring with a unity element 1 having the characteristice 2. Let $M=M_{1, 2}(R)$ and $\Gamma =\left\{\left(\begin{array}{c}n.1\\n.1\end{array}\right): n\in Z, \mbox{n is not divisible by 2}\right\}$.
Then $M$ is a $\Gamma $-ring. Let $N=\{(x, x): x\in R\}\subseteq M$. 
\\Now for all $(x, x)\in N$, $(a, b)\in M$ and $\left(\begin{array}{c}n\\n\end{array}\right)\in \Gamma $, 
we have \begin{eqnarray*}
\lefteqn{(x, x)\left(\begin{array}{c}n\\n\end{array}\right)(a,b)-(a, b)\left(\begin{array}{c}n\\n\end{array}\right)(x,x)}\\&&=(xna-bnx, xnb-anx)\\&&=(xna-2bnx+bnx, bnx-2anx+xna)\\&&=(xna+bnx, bnx+xna)\in N.
\end{eqnarray*}   
Therefore, $N$ is a Lie ideal of $M$.    
\end{ex}

{\bf Acknowledgements:} The research of FH was supported by University Grant Commission of Bangladesh.

\end{document}